\documentclass[12pt,a4paper]{article}
\usepackage{latexsym}
\usepackage{amsfonts}
\usepackage[all]{xy}
\usepackage{amssymb, mathrsfs, amsfonts, amsmath}
\usepackage{amssymb,fontenc,mathrsfs,bbm}
\usepackage{amsbsy}
\newcommand{\noin}{\noindent}

\newcommand{\OM}{{\Omega}}

\newcommand{\POM}{{\partial \, \Omega}}

\newtheorem{theorem}{\bf Theorem}[section]

\newtheorem{definition}[theorem]{\bf Definition}

\begin{document}

\title{A note on the first eigenvalue of spherically symmetric manifolds}
\author{Cleon S. Barroso \and G. Pacelli Bessa\thanks{Research partially supported by CNPq-Brasil}}
\date{\today}
\maketitle
\begin{abstract}
\noin  We give lower and upper bounds for the first eigenvalue of
geodesic balls in spherically symmetric manifolds. These lower and
upper bounds are  $C^{0}$-dependent on the metric coefficients. It
gives better lower bounds for the first eigenvalue of spherical caps
than those from Betz-Camera-Gzyl.

\vspace{.2cm}


\end{abstract}
\section{introduction}Let $B_{\mathbb{N}^{n}(\kappa)}(r) $ be a
 geodesic ball of radius $r>0$
in the  simply connected $n$-dimensional space form
$\mathbb{N}^{n}(\kappa)$
 of constant
sectional curvature $\kappa$ and let
$\lambda_{1}(B_{\mathbb{N}^{n}(\kappa)}(r)$ be its first Laplacian
eigenvalue, i.e. the
 smallest real number $\lambda=\lambda_{1}(B_{\mathbb{N}^{n}(\kappa)}(r))$ for which
 there exists a  function, called a first eigenfunction,   $u\in
 C^{2}(B_{\mathbb{N}^{n}(\kappa)}(r))\cap C^{0}(\overline{B_{\mathbb{N}^{n}(\kappa)}(r)})\setminus\{0\}$,
  satisfying $\triangle u
 +\lambda u=0$ in $B_{\mathbb{N}^{n}(\kappa)}(r)$ with $u\vert \partial
 B_{\mathbb{N}^{n}(\kappa)}(r)=0$. In the case $\kappa=0$, it is well known that
$\lambda_{1}(B_{\mathbb{R}^{n}}(r))=(c(n)/r)^{2}$, where $c(n)$ is
the first zero of the Bessel function  $J_{n/2 -1}$. In the case
$\kappa=-1$, there are fairly good lower and upper bounds for
$\lambda_{1}(B_{\mathbb{H}^{n}}(r))$. For instance, one has  the
upper bounds
$$\begin{array}{lll} \sqrt{\lambda_{1}(B_{\mathbb{H}^{n}}(r))}& \leq & (n-1)(\coth(r/2)-1)/2+   [(n-1)^{2}/4  \,+\\
                              &  & \\
                              &&
+4\pi^{2}/r^{2}+ (n-1)^{2}(\coth(r/2)-1)^{2}/4]^{1/2}.\nonumber
\end{array}$$ See \cite{chavel}, page 49 and see \cite{gage} for sharper upper bounds.
 Moreover, one has the lower bound,
$$\sqrt{\lambda_{1}(B_{\mathbb{H}^{n}}(r))}\geq   \frac{(n-1)\coth(
r)}{2}\cdot $$ See   \cite{yau1}. This lower bound was improved  by
Bessa and Montenegro in \cite{bessa-montenegro} to
\begin{equation}\sqrt{\lambda(r,n)}\geq  \max
\left\{\frac{n}{2 r},\,  \frac{(n-1)\coth( r)}{2}\right\}.
\end{equation} The  case $c=1$  is more delicate. Although the sphere is a well studied manifold,
the  values of  the  first Laplacian eigenvalue
$\lambda_{1}(B_{\mathbb{S}^{n}}(r))$,  (Dirichlet boundary data if
$r<\pi$) are pretty much unknown, with the exceptions
$\lambda_{1}(B_{\mathbb{S}^{n}}(\pi/2))=n$,
$\lambda_{1}(B_{\mathbb{S}^{n}}(\pi))=0$. In dimension two and three
there are good lower bounds due to Barbosa-DoCarmo
\cite{Barbosa-do-Carmo}, Pinsky   \cite{pinsky},  Sato \cite{Sato}
and Friedland-Hayman \cite{friedland-hayman}. In higher dimension,
the lower bounds known (to the best of our knowledge)  are the
following lower bounds due to
 Betz, Camera
and Gzyl  obtained  in \cite{betz-camera-gzyl} via probabilistic
methods.
\begin{equation}\label{eqBCG}\left(\frac{c(n)}{r}\right)^{2}>\lambda_{1}(B_{\mathbb{S}^{n}}(r))\geq
\displaystyle\frac{1}{\displaystyle\int_{0}^{r}\left[\displaystyle\frac{1}{\sin^{n-1}(\sigma)}\cdot
\int_{0}^{\sigma}\sin^{n-1} (s)ds\right]\,d\sigma}\cdot
\end{equation} The upper bound is due to  Cheng's eigenvalue
comparison theorem
 \cite{cheng} since
the Ricci curvature of the sphere is positive (in fact, it needed
only to be
non-negative).\\

In order to state our result let us recall the definition of a
spherically symmetric manifold.
 Let $M$ be a  Riemannian manifold and a point $p\in M$.  For
  each vector
  $\xi \in T_{p}M$, let $\gamma_{\xi}$ be the unique geodesic satisfying
  $\gamma_{\xi}(0)=p$,
  $\gamma_{\xi}'(0)=\xi$ and $d(\xi)=\sup\{t>0:
  {\rm dist}_{M}(p,\gamma_{\xi}(t))=t\}$. Let
  ${\cal D}_{p}=\{t\,\xi\in T_{p}M: 0\le t<d(\xi), \,\vert
   \xi\vert=1\}$ be  the largest open subset of
  $T_{p}M$  such that for
  any $\xi \in {\cal D}_{p}$ the  geodesic
  $\gamma_{\xi}(t)=\exp_{p}(t\,\xi)$ minimizes the distance from
  $p$ to $\gamma_{\xi}(t)$ for all $t\in [0,d(\xi)]$. The cut locus of $p$ is the set ${\rm Cut }(p)=\{ \exp_{p}(d(\xi)\,\xi),\, \xi \in T_{p}M,\, \vert
  \xi \vert =1\}$ and  $M=\exp_{p}({\cal D}_{p})\cup {\rm Cut
  }(p)$.

  The exponential map $\exp_{p}:{\cal D}_{p}\to \exp_{p}({\cal
  D}_{p})$ is a diffeomorphism and is called the geodesic
  coordinates of $M\setminus {\rm Cut}(p)$.
 Fix a vector $\xi \in T_{p}M$, $\vert \xi \vert =1$ and
   denote by $\xi^{\perp}$  the orthogonal
complement of $\{\mathbb{R}\xi\}$ in $T_{p}M$ and let
$\tau_{t}:T_{p}M\to T_{\exp_{p}(t\,\xi)}M$ be the parallel
translation along $\gamma_{\xi}$.  Define the path of linear
transformations $$ {\cal A}(t,\xi):\xi^{\perp}\to\xi^{\perp}$$ by
$${\cal A}(t,\xi)\eta=(\tau_{t})^{-1}Y(t)$$ where $Y(t)$ is the
Jacobi field along $\gamma_{\xi}$ determined by the initial data
$Y(0)=0$, $(\nabla_{\gamma_{\xi}'}Y)(0)=\eta$. Define the map
$${\cal R}(t):\xi^{\perp}\to \xi^{\perp} $$  by $${\cal
R}(t)\eta=(\tau_{t})^{-1}\, {\rm
R}(\gamma_{\xi}'(t),\tau_{t}\,\,\eta)\gamma_{\xi}'(t),$$  where
${\rm R}$ is the Riemann curvature tensor of $M$. It turns out that
the map ${\cal R}(t)$ is a self adjoint map and the path of linear
transformations ${\cal A}(t,\xi)$ satisfies the Jacobi equation
${\cal A}''+{\cal R}{\cal A}=0$ with initial conditions ${\cal
A}(0,\xi)=0$, ${\cal A}'(0,\xi)=I$. On the set $\exp_{p}({\cal
D}_{p})$ the Riemannian metric of $M$ can be expressed by
\begin{equation}\label{eqmetricGeoCoord}ds^{2}(\exp_{p}(t\,\xi))=dt^{2}+\vert {\cal
A}(t,\xi)d\xi\vert^{2}. \end{equation}
\begin{definition}A manifold $M$ is said to be
 spherically symmetric  if the matrix ${\cal
A}(t,\xi)=f(t)I$, for a  function $f\in C^{2}([0, R])$, $R\in (0,
\infty]$ with $f(0)=0$,
 $f'(0)=1$, $f\vert (0, R)>0$.
\end{definition}
 The class of spherically
symmetric manifolds includes the canonical space forms
$\mathbb{R}^{n}$, $\mathbb{S}^{n}(1)$ and $\mathbb{H}^{n}(-1)$. The
$n$-volume $V(r)$ of a geodesic ball $B_{M}(r)$ of radius $r$ in a
spherically symmetric manifold is given by
$V(r)=w_{n}\int_{0}^{r}f^{n-1}(s)ds$, whereas the $(n-1)$-volume
$S(r)$ of the boundary $\partial B_{M}(r)$ is given by
$S(r)=w_{n}f^{n-1}(r)$. Here $w_{n}$ denotes the $(n-1)$-volume of
the sphere $\mathbb{S}^{n-1}(1)\subset \mathbb{R}^{n}$.
 The authors \cite{barroso-bessa} obtained using   fixed point
 methods
the following lower bound for the first eigenvalue
$\lambda_{1}(B_{M}(r))$ of geodesic balls $B_{M}(r)$ with radius $r$
in a spherically symmetric manifold $M$,
\begin{equation}
\lambda_{1}(B_{M}(r))\geq
\frac{1}{\displaystyle\int_{0}^{r}\frac{V(\sigma)}{S(\sigma)}
\,d\sigma}\cdot\label{eqBarroso-Bessa}
\end{equation}It is worth mentioning that this lower bound (\ref{eqBarroso-Bessa}) is
 Betz-Camera-Gzyl's lower bound
 when $M=\mathbb{S}^{n}$. The
purpose of this note is give upper and better lower bounds for
$\lambda_{1}(B_{M}(r))$. We prove the following theorem.
\begin{theorem}\label{thm1}Let $B_{M}(r)\subset M$ be a  ball in
a spherically symmetric Riemannian manifold  with  metric $dt^{2}+
f^{2}(t)d\theta^{2}$, where $f\in C^{2}([0,R])$  with
  $f(0)=0$, $f'(0)=1$, $f(t)>0$ for
all $t\in (0,R]$.
 For every non-negative function  $u\in C^{0}([0,r])$  set $$h(t,u)=\left[u(t)/\int_{t}^{r}\displaystyle\int_{0}^{\sigma}
\left(\frac{f (s) }{f (\sigma)}\right)^{n-1} u(s)dsd\sigma
\right].$$Then
\begin{equation}\label{eqbarroso-bessa}
\sup_{t}h(t,u) \geq \lambda_{1}(B_{M}(r))\geq \inf_{t}h(t,u)
\end{equation}Equality if (\ref{eqbarroso-bessa}) if  and only if $u$
is a first positive
 eigenfunction of
$B_{M}(r)$ and $ \lambda_{1}(B_{M}(r))= h(t,u).$
\end{theorem}In the following table we compare our estimates for
 $\lambda_{1}(r)=\lambda_{1}(B_{\mathbb{S}^{n}}
 (r))$ for $n=2,3$, $r=\pi/8,\pi/4, 3\pi/8, \pi/2, 5\pi/8$
 taking $u(t)=\cos (t\pi/2r)$ with the estimates obtained by
 Betz-Camera-Gzyl.

$$
\begin{tabular}{|l|l|l|l|l|l|}
\hline $n=2/r$ &$\pi/8$&  $\pi/4)$ & $\pi/8$ &$\pi/2$ & $5\pi/8$
\\
\hline BCG/$\lambda_{1}(r)$ & $\geq 25.77 $ & $ \geq6.31 $ & $\geq  2.70$ & $ \geq 1.44$& $ \geq0.85$\\
\hline BB/$\lambda_{1}(r)$& $ \geq35.85$  & $ \geq8.78$ & $ \geq3.76
$ & $ =2 $& $\geq1.01$
\\\hline
\hline $n=3/r$ &$\pi/8$& $\pi/4$ & $3\pi/8$ &$\pi/2$ & $5\pi/8$
\\
\hline BCG/$\lambda_{1}(r) $& $ \geq38.50 $ & $\geq 9.31 $ & $ \geq 3.90$ & $\geq  2$& $\geq 1.10$\\
\hline BB/$\lambda_{1}(r)$ & $\geq 57.94$  & $\geq14.01$ & $
\geq5.86 $ & $ =3 $&$ \geq1.27$\\\hline
\end{tabular}
$$
\section{Proof of Theorem \ref{thm1}}We start recalling the following theorem due to J.
Barta.
\begin{theorem}[Barta, \cite{barta}]\label{barta} Let $\Omega \subset M$ be a
bounded domain  with  piecewise smooth boundary $\POM$ in a
Riemannian manifold. For any $f\in C^{2}(\OM)\cap
C^{0}(\overline{\OM})$ with
    $f\vert \OM>0$ and $f\vert \partial \OM=0$ one has that
\begin{equation}\label{eqBarta1}\sup_{M} ( -\triangle f/f) \geq
\lambda_{1}(\OM ) \geq \inf_{\OM} (-\triangle f/f).
\end{equation}Equality in (\ref{eqBarta1}) holds if and only if   $f$
is a first eigenfunction of $\OM$. The lower bound inequality needs
only that $f\vert \OM>0$.
\end{theorem}

Let $u\in C^{0}([0,r])$, $u\geq 0$. Define a function $T(u)\in
C^{1}([0,r])$ by
$T(u)(t)=\smallint_{t}^{t}\smallint_{0}^{\sigma}(f(s)/f(\sigma))^{n-1}u(s)ds
d\sigma$. Extend $u$ and $Tu$ radially to $B_{M}(r)$
  by
$\tilde{u} (\exp_{p}(t\,\eta))=u(t)$ and $\tilde{T}(u)
(\exp_{p}(t\,\eta))=T(u)(t)$, for $\eta \in \mathbb{S}^{n-1}$.
Observe that $\tilde{T}(u)(\exp_{p}(t\,\eta))\geq 0$, with
$\tilde{T}(u)(\exp_{p}(t\,\eta))= 0$ if and only if $t=r$. We have
that \begin{equation}\label{eqPrincipal}\triangle
\tilde{T}u(\exp_{p}(t\,\eta))=-\tilde{u}(\exp_{p}(t\,\eta))\end{equation}
as a straight forward computation shows. Applying Barta's Theorem we
obtain that
$$ \sup_{t}\frac{u}{T(u)}(t)\geq \lambda_{1}(B(r))\geq
\inf_{t}\frac{u}{T(u)}(t).$$

 Barta's
Theorem says that equality in (\ref{eqPrincipal}) holds if and only
if $\tilde{T}(u)$ is a first eigenfunction. Thus we need only to
show that $\tilde{T}(u)$ is a  first eigenfunction if and only if
$u$ is a  first eigenfunction. Suppose that we have equality in
(\ref{eqPrincipal}) then $\tilde{T}(u)$ is an eigenfunction, this is
 \begin{equation}\label{eqsecundaria}0=\triangle
\tilde{T}u+\lambda_{1}(B_{M}(r))\tilde{T}u=-\tilde{u}+\lambda_{1}(B_{M}(r))\tilde{T}u\end{equation}
Applying the Laplacian in both side of the equation
(\ref{eqsecundaria}) we obtain by equation  (\ref{eqPrincipal}) that
\begin{equation}
0=-\triangle \tilde{u}+\lambda_{1}(B_{M}(r))\triangle
\tilde{T}u=-(\triangle \tilde{u}+\lambda_{1}(B_{M}(r)) \tilde{u})
\end{equation}Therefore $u$ is a first  eigenfunction
with $\lambda_{1}(B_{M}(r))=\displaystyle\frac{u}{T(u)}\cdot$


\begin{thebibliography}{10}
\bibitem{Barbosa-do-Carmo}J. L. Barbosa, M.   Do Carmo,
\newblock{\em
Stability of minimal surfaces and eigenvalues of the Laplacian.}
Math. Z. {\bf 173},  13-28, (1980)
\bibitem{barroso-bessa} C. S. Barroso, G. P. Bessa, {\em  Lower bounds for the first Laplacian eigenvalue of geodesic balls of spherically
symmetric manifolds.} math.DG/0601180. (To appear in Int. J. Appl.
Math. Stat.),


\bibitem{barta}J.  Barta,  {\em Sur la vibration fundamentale
d'une membrane.} C. R. Acad. Sci. {\bf 204},  472-473, (1937).
\bibitem{bessa-montenegro} G. P. Bessa, J.
F. Montenegro. {\em Eigenvalue estimates for submanifolds  with
locally bounded mean  curvature. } \textit{Ann. Global Anal. and
Geom.} \textbf{24}, 279--290, (2003).
\bibitem{betz-camera-gzyl} C. Betz, G. A. Camera, H. Gzyl,     {\em Bounds for the first eigenvalue of a spherical cup. }
 Appl. Math. Optm. {\bf 10}, 193-202, (1983).
\bibitem{chavel}I. Chavel,  {\em Eigenvalues in Riemannian
Geometry.}  Pure and Applied Mathematics,  Academic Press, (1984).
\bibitem{cheng}S. Y. Cheng, {\em Eigenvalue comparison theorems and its
geometric applications.} Math. Z. \textbf{143}, 289--297, (1975).

 \bibitem{friedland-hayman} S. Friedland, W. K. Hayman, {\em
Eigenvalue inequalities for the dirichlet problem on spheres and the
growth of subharmonic functions.} Comment. Math. Helvetici {\bf 51},
133-161, (1976).
 \bibitem{gage} M. Gage, {\em Upper bounds for the first eigenvalue
of the Laplace-Beltrami operator.} \textit{Indiana Univ. Math. J.}
\textbf{29}, 897-912, (1980).


\bibitem{pinsky}  M. A. Pinsky, {\em The first eigenvalue of a sphericall cap} \textit{Appl. Math. Opt.}
\textbf{7},  137-139, (1981).
\bibitem{Sato}S. Sato, {\em Barta's inequalities and the first eigenvalue of a cap domain of a $2$-sphere.}
Math. Z. \textbf{181}, 313-318, (1982).
\bibitem{yau1} S. T. Yau,  {\em Isoperimetric constant and the
first eigenvalue of a compact Riemannian manifold.} Ann. Sci.
\'{E}cole Norm. Sup. $\textbf{4}$  s\'{e}rie t.8,  487-507, (1975).
\end{thebibliography}
\end{document}